\documentclass[12pt]{article}
\usepackage{latexsym,amsfonts,amsmath,graphics,amssymb}
\usepackage{verbatim}
\usepackage{epsfig}

\usepackage[english]{babel}

\usepackage{color,graphicx}
\usepackage{latexsym,amsfonts,amsmath,amssymb,epsfig,amsthm}
\usepackage{euler}

\usepackage[colorlinks=true,linkcolor=black,citecolor=black]{hyperref}

\newtheorem{algorithm}{Algorithm}

\allowdisplaybreaks

\newcommand{\bst}{\boldsymbol{t}}

\newcommand{\bsu}{\boldsymbol{u}}
\newcommand{\bsx}{\boldsymbol{x}}
\newcommand{\bsy}{\boldsymbol{y}}
\newcommand{\bsz}{\boldsymbol{z}}

\newcommand{\bszero}{\boldsymbol{0}}

\newcommand{\bsone}{\boldsymbol{1}}

\newcommand{\rd}{\mathrm{\, d}}

\makeatletter
\newcommand{\rdots}{\mathinner{\mkern1mu\lower-1\p@\vbox{\kern7\p@\hbox{.}}
\mkern2mu \raise4\p@\hbox{.}\mkern2mu\raise7\p@\hbox{.}\mkern1mu}}
\makeatother

\begin{document}

\title{\scshape Applications of geometric discrepancy in numerical analysis and statistics}

\author{Josef Dick\footnote{School of Mathematics and Statistics, The University of New South Wales, Sydney, 2052 NSW, Australia; email: josef.dick@unsw.edu.au} } 

\date{}
\maketitle

\begin{abstract}
In this paper we discuss various connections between geometric discrepancy measures, such as discrepancy with respect to convex sets (and convex sets with smooth boundary in particular), and applications to numerical analysis and statistics, like point distributions on the sphere, the acceptance-rejection algorithm and certain Markov chain Monte Carlo algorithms.
\end{abstract}


{\bf MSC Class:} 65D30, 65D32

\section{Introduction}

The local discrepancy function of a point set $P_{N,s} = \{\bsx_0, \bsx_1, \ldots, \bsx_{N-1}\} \subset [0,1]^s$ measures the difference of the empirical distribution from the uniform distribution with respect to some test sets $\mathcal{A} \subset \mathcal{P}([0,1]^s)$, where $\mathcal{P}([0,1]^s)$ denotes the power set of $[0,1]^s$, that is
\begin{equation*}
\Delta_A(P_{N,s}) = \frac{1}{N} \sum_{n=0}^{N-1} 1_A(\bsx_n) - \lambda_s(A),
\end{equation*} 
where $A \in \mathcal{A}$, $1_{A}$ is the indicator function of the set $A$ and $\lambda_s$ denotes the $s$-dimensional Lebesgue measure. The supremum of $|\Delta_A(P_{N,s})|$ over all sets $A \in \mathcal{A}$ is called the star-discrepancy of $P_{N,s}$ (with respect to the test sets $\mathcal{A}$)
\begin{equation*}
D^\ast_{\mathcal{A}}(P_{N,s}) = \sup_{A \in \mathcal{A}} |\Delta_A(P_{N,s})|.
\end{equation*}

Depending on the choice of test sets, one gets different types of convergent behavior. One well studied example of test sets is that of boxes anchored at the origin $\mathcal{B} = \{[\bszero, \bst): \bst \in [\bszero, \bsone]\}$, where $\bszero = (0, 0, \ldots, 0)$, $\bst = (t_1, t_2, \ldots, t_s) \in [0,1]^s$ and $[\bszero, \bst) = \prod_{j=1}^s [0, t_j)$. In this case upper and lower bounds are known, as well as explicit constructions of point sets which match the best known upper bounds. Variations of anchored boxes, like boxes anchored at different places, boxes which are not anchored or boxes on the torus $(\mathbb{R}/\mathbb{Z})^s$ are all similar and many results are also known in these cases \cite{BC, DP10, DT, KN, Mat}. If the test sets are more complicated, like all convex sets or all convex sets with smooth boundary, then the situation is more complicated. Upper and lower bounds are known, but the upper bounds are often based on probabilistic arguments and are therefore not constructive.

In this paper we provide examples of applications which naturally yield problems in discrepancy theory. The case of anchored boxes is the best understood example of these and provides a connection of low-discrepancy point sets to applications in numerical integration. The connection of discrepancy with respect to other types of test sets is less well-known. We relate point distributions on the sphere, points transformed via inversion to different distributions and spaces, the acceptance-rejection algorithm and Markov chain Monte Carlo algorithm to various discrepancy measures of point sets in the unit cube. This provides a motivation for studying discrepancy with respect to various test sets in the cube.

\section{Numerical integration in the unit cube}

We explain how numerical integration in the unit cube using equal weight quadrature rules leads one to discrepancy with respect to anchored boxes $\mathcal{B}$. A simple explanation in one dimension is the following. Assume that $f$ is absolutely continuous and let $P_{N,1} = \{x_0, x_1, \ldots, x_{N-1} \} \subset [0,1]$. We have
\begin{equation*}
f(x) = f(1) - \int_0^1 1_{[0,t)}(x) f'(t) \rd t
\end{equation*} 
and therefore
\begin{align*}
& \int_0^1 f(x) \rd x - \frac{1}{N} \sum_{n=0}^{N-1} f(x_n) =  \\ & \int_0^1 f'(t) \left[\frac{1}{N} \sum_{n=0}^{N-1} 1_{[0,t)}(x_n)  - \int_0^1 1_{[0,t)}(x) \rd x \right] \rd t = \\ & \int_0^1 f'(t) \Delta_{P_{N,1}}(t) \rd t,
\end{align*} 
where $\Delta_{P_{N,1}}$ is the local discrepancy function given by
\begin{equation*}
\Delta_{P_{N,1} }(t) = \frac{1}{N} \sum_{n=0}^{N-1} 1_{[0,t)}(x_n) - t.
\end{equation*}  
Using H\"older's inequality we therefore get
\begin{equation}\label{Koksma_ineq}
\left| \int_0^1 f(x) \rd x - \frac{1}{N} \sum_{n=0}^{N-1} f(x_n) \right| \le \left(\int_0^1 |f'(t)|^p \rd t)\right)^{1/p} \left(\int_0^1 |\Delta_{P_{N,1} }(t)|^q \rd t\right)^{1/q},
\end{equation}
for H\"older conjugates $1 \le p, q \le \infty$, with the obvious modifications for $p$ or $q = \infty$. Inequality \eqref{Koksma_ineq} is a variation of an inequality due to Koksma~\cite{Koksma}.

From these considerations, one obtains the $L^q$ discrepancy as a quality criterion for the point set $P_{N,1}$:
\begin{equation*}
L^q_{\mathcal{B}}(P_{N,1}) =  \left( \int_{[0,1]} |\Delta_{P_{N,1}}(t)|^q \rd t \right)^{1/q} \quad \mbox{for } 1 \le q \le \infty,
\end{equation*} 
again with the obvious modifications for $q=\infty$.

There is a natural generalization of the above approach to dimensions $s > 1$ by using partial derivatives of $f$. This leads one to discrepancy measures with respect to anchored boxes. Let a point set $P_{N,s} = \{\bsx_0, \bsx_1,\ldots, \bsx_{N-1} \} \subset [0,1]^s$ be given and let  $\bst = (t_1,\ldots, t_s) \in [0,1]^s$. Then we define the local discrepancy function by $$\Delta_{P_{N,s}}(\bst) = \frac{1}{N} \sum_{n=0}^{N-1} 1_{[\bszero,\bst]}(\bsx_n) - \prod_{j=1}^s t_j,$$  where $[\bszero,\bst] = \prod_{j=1}^s [0,t_j]$. Again, by taking the $L^q$ norm of the local discrepancy function, we obtain the $L^q$ discrepancy with respect to anchored boxes $\mathcal{B}$ given by 
\begin{equation*}
L^q_{\mathcal{B}}(P_{N,s}) =  \left( \int_{[0,1]^s} |\Delta_{P_{N,s}}(\bst)|^q \rd \bst \right)^{1/q} \quad \mbox{for } 1 \le q \le \infty,
\end{equation*}
with obvious modifications for $q=\infty$.

There is a generalization of \eqref{Koksma_ineq}, which is due to Hlawka~\cite{Hlawka}. Let $f: [0,1]^s \to \mathbb{R}$ and $P_{N,s} = \{\bsx_0, \bsx_1, \ldots, \bsx_{N-1} \} \subset [0,1]^s$, then 
\begin{align}\label{Hlawka_ineq}
\left| \int_{[0,1]^s} f(\bsx) \rd \bsx - \frac{1}{N} \sum_{n=0}^{N-1} f(\bsx_n) \right| \le & \|f\|_p L^q_{\mathcal{B}}(P_{N,s}),
\end{align} 
where $\frac{\partial^{|u|} f}{\partial \bst_u}(\bst_{u}, \bsone) = 0$ for all $u \subsetneq \{1, 2, \ldots, s\}$ and  
\begin{equation*}
\|f\|_p^p =  \int_{[0,1]^s}  \left| \frac{\partial^s f}{\partial \bst}(\bst) \right|^p \rd \bst.
\end{equation*} 

Several important variations of \eqref{Hlawka_ineq} are known, see for instance Hickernell~\cite{Hi98} and Sloan and Wo\'zniakowski~\cite{SW98} (but these are not discussed here in further detail).

The $L^q$ discrepancy has been intensively studied and many precise results are known. Lower bounds by Roth~\cite{Roth} and Schmidt~\cite{Schmidt} and upper bounds via explicit constructions by Chen and Skriganov~\cite{CS02} and Skriganov~\cite{Skri06} show that
\begin{equation*}
L^q_{\mathcal{B}}(P_{N,s}) \asymp \frac{(\log N)^{\frac{s-1}{2}}}{N} \quad \mbox{for } 1 < q < \infty.
\end{equation*}
The endpoint cases $q=1$ and $q=\infty$ are still open. The following lower bounds are due to Hal\'asz~\cite{Halasz} for $q=1$ and Bilyk and Lacey~\cite{BL08} and Bilyk, Lacey and Vagharshakyan~\cite{BLV08} for $q=\infty$:
\begin{align*}
L^\infty_{\mathcal{B}}(P_{N,s}) \gg_s \frac{(\log N)^{\frac{s-1}{2} + \eta}}{N}, \\ L^1_{\mathcal{B}}(P_{N,s}) \gg_s \frac{\log N}{N}.
\end{align*}

Explicit constructions of point sets are known in each case, see Halton~\cite{Halton}, Hammersley~\cite{Ham}, Sobol~\cite{Sobol}, Faure~\cite{Faure}, Niederreiter~\cite{nie88}, Niederreiter-Xing~\cite{NX95, NX96, NXbook}, Chen-Skriganov~\cite{CS02}, Skriganov~\cite{Skri06}, D.-Pillichshammer~\cite{DP13}, D.~\cite{D13} and others, which show that
\begin{align*}
L^\infty_{\mathcal{B}}(P_{N,s}) \ll_s \frac{(\log N)^{s-1} }{N}, \\ L^1_{\mathcal{B}}(P_{N,s}) \ll_s \frac{(\log N)^{\frac{s-1}{2}} }{N}.
\end{align*}

In the following we consider generalizations of discrepancy measures, which we then relate to problems from numerical analysis and statistics.

\subsection*{Generalizations of the discrepancy with respect to anchored boxes}

The $L^\infty$ discrepancy is often called the star-discrepancy and is denoted by $D^\ast$. Consider the star-discrepancy with respect to anchored boxes $\mathcal{B}$
\begin{equation*}
D^\ast_{\mathcal{B}}(P_{N,s}) = \sup_{\boldsymbol{t} \in [0,1]^s} \left| \frac{1}{N} \sum_{n=0}^{N-1} 1_{[\boldsymbol{0}, \boldsymbol{t})}(\boldsymbol{x}_n) - \prod_{j=1}^s t_j \right|.
\end{equation*}
The sepremum over the boxes $[\bszero, \bst)$ can be replaced by other test sets. This will yield other discrepancy criteria. For instance, the isotropic discrepancy is defined with respect to convex sets $\mathcal{C}$ of $[0,1]^s$. The local isotropic discrepancy is in this case defined by $$\Delta_{P_{N,s}}(C) = \frac{1}{N} \sum_{n=0}^{N-1} 1_{C}(\bsx_n) - \lambda_s(C),$$ where $C \in \mathcal{C}$ is a convex set and $\lambda_s$ is the $s$-dimensional Lebesgue measure. The  isotropic discrepancy is then defined by $$D^\ast_{\mathcal{C}}(P_{N,s}) = \sup_{C \in \mathcal{C}} \left|\Delta_{P_{N,s}}(C) \right|.$$ The connection to numerical integration is not as clear in this case as for the case of anchored boxes.

Again, a number of results are known about the isotropic discrepancy due to  Beck~\cite{Be88}, Hlawka~\cite{Hlawka},   Laczkovich~\cite{La95}, M\"uck and Philipp~\cite{MP75}, Niederreiter~\cite{nie72a, nie72b}, Niederreiter and Wills~\cite{NW75}, Schmidt~\cite{Sch75}, Stute~\cite{Stute}, Zaremba~\cite{Za70}. In dimension $2$ it is known that 
\begin{align*}
N^{-\frac{2}{3}} \ll D^\ast_{\mathcal{C}}(P_{N,2}) \ll N^{-\frac{2}{3}} (\log N)^4,
\end{align*}
whereas in dimension $s > 2$ we have
\begin{align*}
N^{-\frac{2}{s+1}} \ll_s D^\ast_{\mathcal{C}}(P_{N,s}) \le (D^\ast(P_{N,s}))^{1/s}.
\end{align*}
Using known constructions of point sets with small star-discrepancy, one obtains the upper bound
\begin{equation*}
D^\ast_{\mathcal{C}}(P_{N,s}) \le C_s \frac{\log N}{N^{1/s}}.
\end{equation*}
The construction of the point sets is explicit in this case. However, there is a gap between the upper and lower bound and the precise rate of convergence remains unknown.

In the next section we consider numerical integration over the unit sphere.

\section{Numerical integration over the unit sphere}

A natural way to define test sets that takes the symmetry of the unit sphere $\mathbb{S}^s =  \{ \boldsymbol{x} \in \mathbb{R}^{s+1}: \|\boldsymbol{x}\| = 1\}$ into account is to use spherical caps. Spherical caps with center $\bsx$ and height $t$ are defined by
\begin{equation*}
C(\bsx, t) = \{\bsz \in \mathbb{S}^s: \langle \bsz, \bsx \rangle > t\}.
\end{equation*}
Let $\mathcal{S} = \{C(\bsx,t): \bsx \in \mathbb{S}^s, -1 \le t \le 1\}$ denote the set of spherical caps.
The spherical cap $L^q$ discrepancy is now given by
\begin{equation*}
L_{\mathcal{S}}^q(P_{N,s}) = \left(\int_{-1}^1 \int_{\mathbb{S}^s} \left| \frac{1}{N} \sum_{n=0}^{N-1} 1_{C(\bsx,t)}(\bsx_n) - \sigma_s(C(\bsx,t)) \right| \rd \sigma_s(\bsz) \rd t \right)^{1/q},
\end{equation*}
where  $\sigma_s$ is the normalized surface Lebesgue measure on the sphere $\mathbb{S}^s$. As in the cube case, the spherical cap discrepancy is related to numerical integration of functions on the sphere. The Koksma-Hlawka type inequality is of the form 
\begin{equation*}
\left| \int_{\mathbb{S}_s} f(\bsx) \rd \sigma_s(\bsx) - \frac{1}{N} \sum_{n=0}^{N-1} f(\bsx_n) \right| \le L^q_{\mathcal{S}}(P_{N,s}) \|f\|_p,
\end{equation*}
where $\|f\|_p$ is a suitable function norm, see Brauchart and D.~\cite{BD13}.

Bounds on the spherical cap discrepancy have been established by Beck~\cite{Be84},  Schmidt~\cite{Sch69} and Stolarsky~\cite{Sto73}. It is known that
\begin{equation*}
N^{-\frac{1}{2} - \frac{1}{2s}} \ll_s L_{\mathcal{S}}^q(P_{N,s}) \ll_s \sqrt{\log N} N^{-\frac{1}{2} - \frac{1}{2s}}.
\end{equation*}
However, the upper bound is based on probabilistic arguments and there are no known optimal explicit constructions satisfying this upper bound.

A number of non-optimal results follow from the work of Grabner and Tichy~\cite{GT93}, Lubotzky, Philipps and Sarnak~\cite{LPS1, LPS2} and Aistleitner, Brauchart and D.~\cite{ABD12}.

We briefly discuss the explicit constructions of points on the sphere $\mathbb{S}^2$ from Aistleitner, Brauchart and D.~\cite{ABD12}. The idea there is to use a transformation from the square $[0,1]^2$ to the sphere $\mathbb{S}^2$ which preserves the measure. The so-called Lambert transform  $\Phi: [0,1]^2 \to \mathbb{S}^2$ given by
\begin{equation*}
\Phi(x,y) = \left( 2 \cos(2\pi x) \sqrt{y-y^2}, 2 \sin(2 \pi x) \sqrt{y - y^2}, 1-2y  \right)
\end{equation*}
has this property, i.e., for any Lebesgue measurable set $J \subseteq [0,1]^2$ we have $\lambda_2(J) = \sigma_2(\Phi(J))$, where $\Phi(J) = \{\Phi(\bsx): \bsx \in J\}$.

\begin{figure}[ht]
\begin{center}
\includegraphics[scale=0.33]{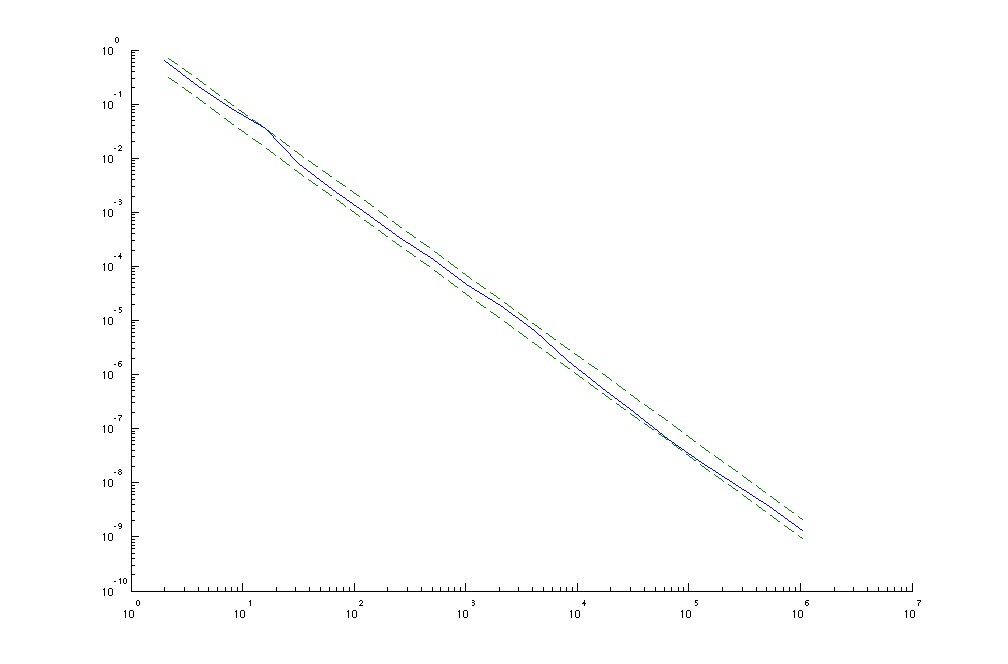}
\caption{\label{fig1} The dashed lines show $N^{-3/2}$ and $(9/4) N^{-3/2}$, and the curve shows the squared spherical cap $L^2$ discrepancy, where the quadrature points are a digital net mapped to the sphere.}
\end{center}
\end{figure}

In order to obtain points on the sphere $\mathbb{S}^2$, we proceed in the following way. We map the points $\{\bsx_0, \bsx_1, \ldots, \bsx_{N-1} \} \in [0,1]^2$ to $\Phi(\bsx_0), \Phi(\bsx_1), \ldots, \Phi(\bsx_{N-1}) \in \mathbb{S}^2$. It was shown in Aistleitner, Brauchart and D.~\cite{ABD12} that if $\bsx_0, \bsx_1, \ldots, \bsx_{N-1}$ have low `discrepancy' with respect to anchored boxes in the square, then $\Phi(\bsx_0), \Phi(\bsx_1), \ldots, \Phi(\bsx_{N-1})$ have low spherical cap discrepancy. Figure~\ref{fig1} shows some numerical results of a digital net mapped to the sphere $\mathbb{S}^2$. The result indicates that these point sets achieve the optimal rate of convergence of the spherical cap $L^2$ discrepancy.

First note that the optimal rate of convergence for boxes in the square differs from the optimal rate of convergence for spherical caps on the sphere. This is not surprising when considering the inverse sets of spherical caps $B(\bsx,t) = \Phi^{-1}(C(\bsx,t)) = \{\bsz \in [0,1]^2: \Phi(\bsz) \in C(\bsx,t)\}$. These sets have particular shapes which change as $\bsx$ and $t$ vary. They are not convex, however, they can be broken up into a small number of convex parts and parts whose complement with respect to some rectangle is convex. Their boundary is smooth except for the pole caps where $\bsx$ is either the north pole $(1,0,0)$ or south pole $(-1,0,0)$, in which case $B(\bsx,t)$ is a rectangle. The curvature of the boundary is unbounded, which can be seen when $\bsx$ moves to one of the poles where the smooth boundary curve of $B(\bsx,t)$ turns into a rectangle. Thus the sets $B(\bsx,t)$ do not have any discernible features. However, for most sets $B(\bsx,t)$, the boundary is smooth and has bounded curvature. Thus, for the most part, the sets $B(\bsx, t)$ can be described by convex sets with smooth boundary which has bounded curvature. 

In the following we briefly discuss discrepancy in the cube with respect to convex test sets with smooth boundary, since the problem of the discrepancy of points mapped to the sphere using the Lambert transform is, to a large degree, related to this discrepancy.

\subsection*{Discrepancy with respect to convex sets with smooth boundary}

We consider now the star-discrepancy in the torus $(\mathbb{R}\setminus \mathbb{Z})^2$ with respect to the convex test sets $\mathcal{E}$, whose boundary curve is twice continuously differentiable with minimal curvature divided by maximal curvature bounded away from 0. It was shown by Beck and Chen~\cite{BC} that
\begin{equation*}
\frac{1}{\sqrt{\log N} N^{3/4}}\ll L^{\infty}_{\mathcal{E}}(P_{N,2}) \ll \frac{\sqrt{\log N}}{N^{3/4}}
\end{equation*}
A generalization to arbitrary dimension can be found in Drmota~\cite{Drmota}, which shows that
\begin{equation*}
\frac{1}{N^{\frac{1}{2} + \frac{1}{2s} }}\ll_s L^{\infty}_{\mathcal{E}}(P_{N,s}) \ll_s \frac{\sqrt{\log N}}{N^{\frac{1}{2} + \frac{1}{2s} }}.
\end{equation*}
The discrepancy bounds for this case are very similar to the discrepancy bounds for the spherical cap discrepancy on $\mathbb{S}^s$. As in the sphere case, no explicit constructions of point sets achieving the upper bound on the discrepancy with respect to convex sets with smooth boundary is known. Indeed, a solution of one of these problems may also yield a solution to the other problem. The numerical results for the spherical cap $L^2$ discrepancy may indicate that classical low-discrepancy constructions such as digital nets and Fibonacci lattices are optimal. Hence the question arises whether this is also true for the discrepancy in the square with respect to convex sets with smooth boundary, as studied by Beck and Chen~\cite{BC}.

\section{Inverse transformation and test sets}\label{sec_inverse}

Assume now we want to approximate the integral
\begin{equation*}
\int_{G} f(\bsx) \psi(\bsx) \rd \bsx,
\end{equation*}
where $\psi$ is a probability density function on $G \subseteq \mathbb{R}^s$. In some cases there is a mapping $\Phi: [0,1]^s \to \mathbb{G}$ which is measure preserving in the following sense: For every Lebesgue measurable set $A \subseteq [0,1]^s$ we have
\begin{equation*}
\lambda_s(A) = \int_{\Phi(A)} \psi(\bsx) \rd \bsx.
\end{equation*}
See the inverse Rosenblatt transform~\cite{Ro52}, or the monographs by Devroye~\cite{De86} and H\"ormann, Leydold and Derflinger~\cite{HLD04}. In the previous section we saw an example of this situation, however such problems come up in other contexts as well (see for instance Kuo, Dunsmuir, Sloan, Wand, and Womersley~\cite{KDSWW} for an example in the context of quasi-Monte Carlo integration). Now assume that we want to study discrepancy with respect to boxes in $G$. Let $P = \{\bsx_0, \bsx_1, \ldots, \bsx_{N-1}\} \subset [0,1]^s$. Then we can define the discrepancy
\begin{equation*}
\sup_{\bst \in \mathbb{R}^s} \left|\frac{1}{N} \sum_{n=0}^{N-1} 1_{G \cap [-\boldsymbol{\infty}, \bst]}(\Phi(\bsx_n)) - \int_{G \cap [-\boldsymbol{\infty}, \bst]} \psi(\bsx) \rd \bsx \right|.
\end{equation*}
This discrepancy can be translated to a discrepancy in the unit cube
\begin{equation*}
\sup_{A \in \mathcal{A} } \left| \frac{1}{N} \sum_{n=0}^{N-1} 1_A(\bsx_n) - \lambda_s(A) \right|,
\end{equation*}
where $\mathcal{A}$ consists of all sets $A_{\bst} = \{\bsx \in [0,1]^s: \Phi(\bsx) \in G \cap [\bszero, \bst)\}$ for all $\bst \in [0,1]^s$. In general, the sets $A_{\bst}$ are not boxes anymore and therefore one wants to have point sets $\bsx_0, \bsx_1, \ldots, \bsx_{N-1} \in [0,1]^s$ which have small discrepancy with respect to the test sets $\mathcal{A}$ rather than boxes. Such a problem also comes up for instance in L'Ecuyer, Lecot and Tuffin~\cite{KLT} and their arrayRQMC method. In some cases, one can do stratified sampling, that is, divide the cube into subcubes with side length $2^{-k}$ and randomly place a point in each box. Then one can show a convergence rate of order $N^{-1/2 - 1/(2s)}$. However, for instance, in L'Ecuyer, Lecot and Tuffin~\cite{KLT} better rates of convergence where observed when using low-discrepancy point sets. In Kuo, Dunsmuir, Sloan, Wand and Womersley~\cite{KDSWW} it was observed that the choice of transformation $\Phi$ influences the rate of convergence, but it is a priori not clear what choice of $\Phi$ yields the best results. For these types of applications it would be interesting to have point sets which achieve good convergence rates for various types of test sets.

In this context, we mention one construction of explicit point sets in dimension $2$ which considers more general test sets. Namely the construction by Bilyk, Ma, Pipher, Spencer~\cite{BMPS} where discrepancy with respect to certain rotated boxes is considered.

\section{Acceptance-rejection sampler}

In statistical sampling one often wants to sample from a target distribution. The standard procedure to obtain samples from a given distribution is to invert the cumulative distribution function (cdf), which can be used to map points from the cube $[0,1]^s$ to the required domain. However, this is often not possible (or difficult), in which case one has to resort to other methods. As an example, consider the unnormalized density function $\psi(x) = x^2 + \sin (\pi x)$ for $x \in [0,2]$. The normalization constant is $\int_0^2 (x^2 + \sin(\pi x) ) \,\mathrm{d} x= 8/3$. The cdf is given by
\begin{equation*}
\Psi(t) = \frac{3}{8} \int_0^t \psi(x) \,\mathrm{d} x =  \frac{3}{8} \int_0^t (x^2 + \sin(\pi x) ) \,\mathrm{d} x = \frac{t^3}{8} + \frac{1 - \cos (\pi t) }{\pi}.
\end{equation*}
In order to be able to sample directly from the density $3\psi/8$, one would have to invert $\Psi$. Since this is not feasible, one has to resort to other methods like the acceptance-rejection algorithm.

The acceptance rejection sampler proceeds in the following way.
\begin{algorithm}\label{AR_alg}
Let $\psi:[0,1]^s \to \mathbb{R}$ be a (unnormalized) probability density function.
\begin{enumerate}
\item Choose $L > 0$ such that $\frac{\psi(\bsx)}{L} \le 1$ for all $\bsx \in [0,1]^s$.
\item Generate a point set $P_{M,s+1} = \{\bsx_0, \bsx_1, \ldots, \bsx_{M-1}\} \subset [0,1]^{s+1}$. Assume that $\bsx_n = (x_{n,1}, x_{n,2}, \ldots, x_{n,s+1})$.
\item Choose $$I = \{0 \le n < M: \psi(x_{n,1}, x_{n,2}, \ldots, x_{n,s}) \le L x_{n,s+1}\}.$$
\item Return the point set $Q = \{(x_{n,1}, x_{n,2},\ldots, x_{n,s}): n \in I\}$.
\end{enumerate}
\end{algorithm}
It is well known that if $P_{M, s+1} = \{\bsx_0, \bsx_1,\ldots, \bsx_{M-1}\}$ is chosen i.i.d. uniformly distributed, then the point set $Q$ has distribution with law $\widetilde{\psi}$ (where $\widetilde{\psi} = \psi \left(\int_{[0,1]^s} \psi(\bsx) \,\mathrm{d} \bsx \right)^{-1}$ is the normalized density function). For a proof see for instance Robert and Casella~\cite[Section~2.3]{RC04}.

\begin{figure}[ht]
\begin{center}
\includegraphics[scale=0.65]{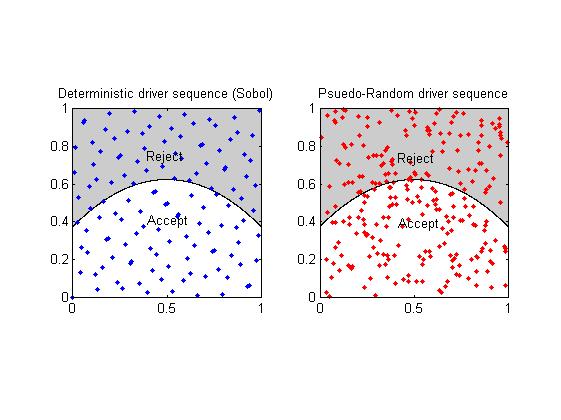}
\end{center}
\caption{\label{fig2} Illustration of the acceptance-rejection algorithm. Points below the curved line are accepted and then projected onto the $x$-axis. The left-hand side uses a deterministic point set $P_{M, s+1}$, whereas the right-hand side uses random samples $P_{M, s+1}$.}
\end{figure}

Numerical tests have been performed by Morokoff and Caflisch~\cite{MC95}, Moskowitz and Caflisch~\cite{MC96} and Wang~\cite{Wang}, where the random point set $P_{M,s+1}$ in Algorithm~\ref{AR_alg} is replaced by a low discrepancy point set with the intention to obtain samples which have better distribution properties. The difference between random point sets and deterministic point sets in the acceptance-rejection algorithm is illustrated in Figure~\ref{fig2}.

Let $Q = \{\bsy_0, \bsy_1, \ldots, \bsy_{N-1}\}$ be the set of points generated by Algorithm~\ref{AR_alg}. To study the performance of this type of algorithm, we introduce the discrepancy
\begin{equation*}
D^\ast(Q) = \sup_{\bst \in [0,1]^{s} } \left| \frac{1}{N}  \sum_{n=0}^{N-1} 1_{[\bszero, \bst)}(\bsy_n) - \frac{1}{C} \int_{[\bszero, \bst)} \psi(\bsx) \rd \bsx \right|,
\end{equation*}
where $C = \int_{[\bszero, \bst]} \psi(\bsx) \rd \bsx$.

\begin{figure}[ht]
\begin{center}
\includegraphics[scale=0.45]{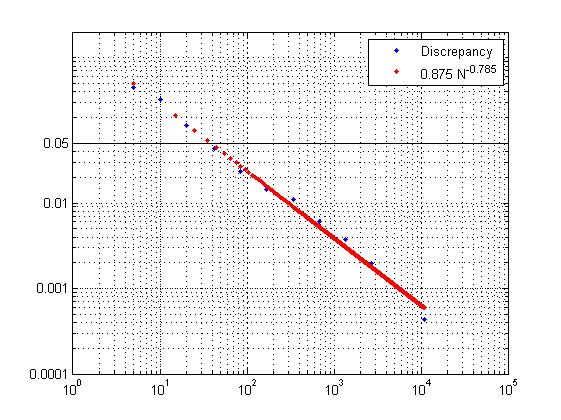}
\end{center}
\caption{\label{fig3} Numerical result of the acceptance-rejection algorithm using low-discrepancy point sets. The convergence rate is approximately of order $N^{-0.7}$, which is better than the rate one would expect when using random samples (which is $N^{-0.5}$).}
\end{figure}

Some simple numerical tests confirm that low-discrepancy point sets in Algorithm~\ref{AR_alg} can improve the performance of the acceptance-rejection algorithm. For instance, in Zhu and D.~\cite{ZD13} the following example was considered: let a unnormalized target density be given by
\begin{equation*}
\psi(x)=\frac{3}{4}-(x-\frac{1}{2})^2, ~~ x\in [0,1].
\end{equation*} 
Figure~\ref{fig3} shows the discrepancy of the point set when the proposal points are a digital net.

The discrepancy $D^\ast(Q)$ can be written in terms of the discrepancy of the proposal points $P_{M,s+1}$. To do so, let
\begin{equation*}
A = \{\bsx = (x_1, x_2,\ldots, x_{s+1}) \in [0,1]^{s+1}: f(x_1, x_2, \ldots, x_s) \le L x_{s+1}\}.
\end{equation*}
Then we have
\begin{equation*}
D^\ast(Q) = \sup_{\bst \in [0,1]^{s+1}} \left| \frac{1}{M} \sum_{n=0}^{M-1} 1_{A \cap [\bszero, \bst]}(\bsx_n) - \lambda_{s+1}(A \cap [\bszero, \bst]) \right|.
\end{equation*}
Thus the discrepancy of the samples obtained from Algorithm~\ref{AR_alg} coincides with the discrepancy of the points $P_{M,s+1}$ with respect to the test sets $A \cap [\bszero, \bst)$. If $\psi$ is smooth and concave, then the  set $A \cap [\bszero, \bst]$ is convex and has smooth boundary, except at the intersection points of the boundaries of $A$ and $[\bszero, \bst]$ and the intersection of the faces of the box. Thus bounds on the discrepancy with respect to convex sets with smooth boundary may also be of help in this problem.

Some results about the discrepancy $D^\ast(Q)$ are known from Zhu and D.~\cite{ZD13}. If $P_{M, s+1}$ is a low-discrepancy point set, then for any concave unnormalized density $\psi$ we have $D^\ast(Q) \ll N^{- \frac{1}{s+1} }$. On the other hand, for every point set $P_{M,s+1}$, there is a concave function $\psi$ such that $D^\ast(Q) \gg N^{-\frac{2}{s+2} }$. Thus, in general, the upper bound cannot be significantly improved without using further assumptions. Smoothness of $\psi$ would be an assumption where one would hope to be able to get better results.

\section{Markov chain Monte Carlo and completely uniformly distributed sequences}

Let $G \subseteq \mathbb{R}^d$ be a state space and $\varphi: G \times [0,1]^s \to G$. Then  one obtains a Markov chain by choosing a starting point $\bsx_0 \in G$, generating a sequence of random numbers $\bsu_0, \bsu_1, \ldots \in [0,1]^s$ and setting $\bsx_n = \varphi(\bsx_{n-1}; \bsu_{n-1})$. In Chen, D. and Owen~\cite{CDO11} the authors studied Markov chains and conditions under which the Markov chain consistently samples a target distribution $\psi$, that is, for every continuous function defined on $G$ we have
\begin{equation}\label{eq_consistent}
\lim_{N\to \infty} \frac{1}{N} \sum_{n=0}^{N-1} f(\bsx_n) = \int_G f(\bsx) \rd \bsx.
\end{equation}
In Chen, D. and Owen~\cite{CDO11} it was shown that if the random numbers $\bsu_0, \bsu_1, \ldots$ are completely uniformly distributed (and some further assumptions on the update function are satisfied), then \eqref{eq_consistent} holds. Complete uniform distribution is a condition on a sequence of numbers $u_0, u_1, u_2, \ldots \in [0,1]$ which ensures statistical independence of successive terms in some sense. Its definition is based on discrepancy and works as follows. Let $\bsu^{(s)}_0 = (u_0, u_1, \ldots, u_{s-1})$, $\bsu^{(s)}_1 = (u_s, u_{s+1}, \ldots, u_{2s-1})$, and so on. Let $\mathcal{U}_{N,s} = \{\bsu^{(s)}_0, \bsu^{(s)}_1, \ldots, \bsu^{(s)}_{N-1}\}$. Then the sequence $u_0, u_1, u_2, \ldots$ is completely uniformly distributed (CUD) if for all dimensions $s \ge 1$ we have $$\lim_{N\to \infty} D_{\mathcal{B}}^\ast(\mathcal{U}_{N,s}) = 0.$$

Explicit constructions of sequences which are completely uniformly distributed have been established by Levin~\cite{Levin} and Shparlinski~\cite{Shparlinski}. For instance, Levin~\cite{Levin} showed error bounds of the form $D_{\mathcal{B}}^\ast(P_{N,s}) \ll N^{-1} (\log N)^{s+\varepsilon}$. However, in the application arising in Chen, D. and Owen~\cite{CDO11}, one needs $s \asymp \log N$. In this case we have $$N^{-1} (\log N)^{s + \varepsilon} \approx N^{-1} (\log N)^{\log N} = N^{-1 + \log \log N}.$$ Thus these results do not guarantee convergence of the discrepancy as $N$ tends to $\infty$. Instead, a different approach is required. In Chen, D. and Owen~\cite{CDO11} the existence of a sequence of numbers was shown for which for all $s \ge 1$ we have
\begin{equation*}
D^\ast_{\mathcal{B}}(\mathcal{U}_{N,s}) \le C \sqrt{s\frac{\log N }{N} }.
\end{equation*}
For $s \asymp \log N$ one therefore gets that $D^\ast_{\mathcal{B}}(\mathcal{U}_{N,s}) \le C \frac{\log N }{\sqrt{N} }$. In Aistleitner and Weimar~\cite{AW13} an improvement was obtained where $D^\ast_{\mathcal{B}}(\mathcal{U}_{N,s}) \le C \sqrt{s\frac{\log \log N }{N} }$. However, no explicit construction of such a sequence is known. Such sequences could be of use in applications in Markov chain Monte Carlo.

\section{Uniformly ergodic Markov chains and push-back discrepancy}

The result in Chen, D. and Owen~\cite{CDO11} yields consistency for certain Markov chains for completely uniformly distributed sequences, but does not yield any convergence rates. This question was addressed in D., Rudolf and Zhu~\cite{DRZ13}. Therein, the discrepancy of the sample points in the state space $G$ with respect to the target distribution was related to the discrepancy of the driver sequence $u_0, u_1, u_2, \ldots \in [0,1]$ with respect to the uniform measure. Again, through the update function, the test sets defined in $G$ are distorted in the unit cube, as in Section~\ref{sec_inverse}. In general, one therefore does not have boxes as test sets anymore. Additionally, one also needs the statistical independence of the driver sequence measured in terms of complete uniform distribution. This yields a generalized definition of completely uniformly distributed point sets, where one does not have boxes as test sets. We call the underlying discrepancy a 'push-back discrepancy', since it is derived from the discrepancy in the state space by inverting the update function. 
We provide some details in the following.

Let $\;\mathcal{U}_{N,s} = \{\bsu_0, \bsu_1,\ldots, \bsu_{N-1} \}
\subset [0,1]^{s}$ denote the points which drive the Markov chain via the update function. Let $\varphi: G \times [0,1]^s \to G$ denote again the update function, so that $\bsx_n = \varphi(\bsx_{n-1}; \bsu_{n-1})$. The $n$ times iterated update function is denoted by $\varphi_n$, that is, we have $\bsx_n = \varphi_n(\bsx_0; \bsu_0, \bsu_1, \ldots, \bsu_{n-1})$. Then we define the sets
\begin{align*}
C_{n,\bsx_0}(A) = & \{ \bsz \in [0,1]^{n s}:  \varphi_n(\bsx_0; \bsz) \in A \},
\end{align*}
for $A\in\mathcal{B}(G)$, the Borel $\sigma$ algebra of $G$, and $n\in \mathbb{N}$. The local discrepancy function of the driver point set is then given by
\[
\Delta_{N,A,\psi,\varphi}(\mathcal{U}_{N,s})
=\frac{1}{N}  \sum_{n=1}^{N} \left[1_{(\bsu_0,\ldots, \bsu_{n-1}) \in C_{n,\psi}(A)} -
\lambda_{n s}(C_{n,\psi}(A)) \right],
\]
and the discrepancy of the driver sequence is given by
\begin{equation*}
D^\ast_{\mathscr{A},\psi,\varphi}(\mathcal{U}_{N,s})
 = \sup_{A \in \mathscr{A}} \left|\Delta_{N,A,\psi,\varphi}(\mathcal{U}_{N,s}) \right|.
\end{equation*}
We call $D^\ast_{\mathscr{A},\psi,\varphi}(\mathcal{U}_{N,s})$ push-back discrepancy of $\mathcal{U}_{N,s}$.

The push-back discrepancy combines two principles: The discrepancy in the cube with respect to general test sets and the principle of complete uniform distribution. As for the discrepancy with respect to general test sets, there are no known explicit constructions of point sets with small push-back discrepancy. However, in D., Rudolf and Zhu~\cite{DRZ13} it was shown that there exist points $\bsu_0, \bsu_1, \ldots, \bsu_{N-1} \in [0,1]^s$ such that
\begin{equation*}
D^\ast_{\mathscr{A}, \psi,\varphi}(\mathcal{U}_n) \ll \sqrt{\frac{\log N}{N}}.
\end{equation*}
Such a point set would have direct applications in Markov chain Monte Carlo.

\subsection*{Acknowledgment}

The author is supported by an ARC Queen Elizabeth II Fellowship and an ARC Discovery Project. The help of Houying Zhu is gratefully acknowledged.

\end{document}